\documentclass[twoside,11pt]{amsart}

\usepackage{amsmath,latexsym,amssymb,mathptm, times}
\usepackage{amsmath,latexsym,amssymb,times,enumerate,xypic, stmaryrd}
\input amssym.def
\input amssym
\input xypic
\input xy
\xyoption{all}
\def\demo{\noindent{\bf Proof. }}
\def\sqr#1#2{{\vcenter{\hrule height.#2pt
        \hbox{\vrule width.#2pt height#1pt \kern#1pt
                \vrule width.#2pt}
        \hrule height.#2pt}}}
\def\square{\mathchoice\sqr64\sqr64\sqr{4}3\sqr{3}3}
\def\QED{\hfill$\square$}

\def\tratto{\mbox{\rule{2mm}{.2mm}$\;\!$}}

\def\m{{\mathfrak m}}

\def\p{{\mathfrak p}}

\newtheorem{Theorem}{Theorem}[section]
\newtheorem{Lemma}[Theorem]{Lemma}
\newtheorem{Corollary}[Theorem]{Corollary}
\newtheorem{Proposition}[Theorem]{Proposition}

\newtheorem{Notation and Discussion}[Theorem]{Notation and Discussion}
\newtheorem{Assumptions and Discussion}[Theorem]{Assumptions and Discussion}

\newtheorem{Example}[Theorem]{Example}

\begin{document}

\baselineskip=16pt

\title[Generalized Hilbert coefficients and  Northcott's inequality]
{\Large\bf Generalized Hilbert coefficients and  Northcott's inequality}

\author[ Y. Xie]
{ Yu Xie}

\thanks{AMS 2010 {\em Mathematics Subject Classification}.
Primary 13D40; Secondary 13A30, 13H15,  13C14, 13C15.}


\address{Department of Mathematics and Statistics, Penn State Altoona,
Altoona, Pennsylvania 16601} \email{yzx1@psu.edu}

\vspace{-0.1in}

\begin{abstract}
Let $R$ be a Cohen-Macaulay local ring of dimension $d$ with infinite residue field. Let  $I$ be  an $R$-ideal that has analytic spread $\ell(I)=d$,  $G_d$ condition and the Artin-Nagata property $AN^-_{d-2}$. We provide a formula relating the length $\lambda(I^{n+1}/JI^{n})$ to the difference $P_I(n)-H_I(n)$, where $J$ is a general minimal reduction of $I$, $P_I(n)$ and $H_I(n)$ are the  generalized Hilbert-Samuel polynomial and the generalized Hilbert-Samuel function  in the sense of C. Polini and Y. Xie.  We then use it to establish formulas to compute the higher generalized Hilbert coefficients of $I$. As an application, we extend Northcott's inequality to non $\m$-primary ideals. When equality holds in the generalized Northcott's inequality,
 the ideal $I$ enjoys nice properties. Indeed, in this case, we prove that
 the reduction number of $I$ is at most one and the associated graded ring of $I$ is Cohen-Macaulay.
 We also recover  results of G. Colom$\acute{{\rm e}}$-Nin, C. Polini, B.  Ulrich and Y. Xie on the positivity of the  generalized first Hilbert coefficient $j_1(I)$.
 Our work extends that of S. Huckaba, C. Huneke and A. Ooishi  to ideals that are not necessarily   $\m$-primary.
\end{abstract}

\maketitle

\vspace{-0.25in}

\section{Introduction}


  Multiplicities and Hilbert functions play important role in  commutative algebra and algebraic geometry. 
  It is well-known that multiplicities are widely used to study intersection theory and singularity theory. Besides that, multiplicities and Hilbert functions  reflect various algebraic and geometric properties of an ideal $I$ in a Noetherian local ring $R$. In particular, they provide useful information on the arithmetical properties, like the depth, of the associated graded ring $G$, where   $G$ ($G={\rm gr}_I(R):=\oplus_{n=0}^{\infty}I^n/I^{n+1}$)
is an
algebraic construction whose projective scheme represents the
exceptional fiber of the blowup of a variety along a subvariety.

The classical multiplicities and Hilbert functions (i.e., the Hilbert multiplicity and the Hilbert function) are only defined  for ideals that are primary to the maximal ideal $\m$ of $R$. In order to study properties associated to non $\m$-primary ideals, one has to define generalized multiplicities and generalized Hilbert functions. One of the generalizations of multiplicities of ideals is called the $j$-multiplicity. It was introduced  by R. Achilles and M. Manaresi in 1993 to study improper intersections of two varieties \cite{AM}. In 1999, H. Flenner, L. O' Carroll and W. Vogel defined the generalized Hilbert function using
the 0th local cohomolodgy functor  \cite[Definition 6.1.5]{FOV}.
In 2003, C. Ciuperc${\rm \breve{a}}$ introduced the generalized Hilbert coefficients  via a different approach -- the bigraded ring ${\rm gr}_{\m}(G)$ \cite{C}.  Recently,  C. Polini and Y. Xie  re-conciliated  both approaches and    defined the concepts of  the generalized Hilbert polynomial and the generalized Hilbert coefficients  following the approach of H. Flenner, L. O' Carroll and W. Vogel \cite{PX1}.
One of the  fundamental properties  proved by C. Polini and Y. Xie illustrates  the behavior of the generalized Hilbert
function
under a  hyperplane section
 \cite{PX1}. Indeed, they proved that  the first $d-1$ generalized Hilbert coefficients $j_0(I), \ldots, j_{d-2}(I)$, where $d={\rm dim}\,R$, are preserved under a general hyperplane section. This nice property allows us to study the generalized Hilbert coefficients by reducing to the lower dimensional case.

  The generalized Hilbert coefficients are important invariants of the ideal $I$.  It is well-known that the normalized leading coefficient $j_0(I)$ (i.e., the {\it $j$-multiplicity} of $I$) can be computed using general elements (by \cite[Theorem 3.8]{AM} and \cite[Corollary 2.5]{X}). This number
 was used to prove the refined Bezout's theorem~\cite{FOV},  to detect integral dependence
 of non $\m$-primary ideals (extension of the fundamental theorem of Rees) \cite{FM1},
 and to study the depth of the associated graded ring of an arbitrary ideal (see \cite{PX} and \cite{MX}).
 The next normalized coefficient $j_1(I)$ is called the {\it   generalized first
 Hilbert coefficient} of $I$. In the $\m$-primary case, $j_1(I)=e_1(I)$
 is called the {\it Chern number}  by W. V. Vasconcelos
for its tracking position in distinguishing
Noetherian filtrations with the same Hilbert multiplicity \cite{V}.
The coefficient $e_1(Q)$, where $Q$ is a parameter ideal,  was used to
characterize the Cohen-Macaulay property for large classes of rings \cite{GGHOPV}.
  Moreover,  G. Colom$\acute{{\rm e}}$-Nin, C. Polini, B.  Ulrich and Y. Xie use  $j_1(I)$, where $I$ is an arbitrary ideal, to bound the number of steps  in a process of  normalization of ideals \cite{NPUX}.
Therefore it
is very important to establish properties such as positivity for the higher generalized Hilbert coefficients.

In the case of  $\m$-primary ideals, there are a number of formulas to compute the
  Hilbert coefficients (see  for instance,  \cite{Huneke} and \cite{Huckaba}).
In 1987, C. Huneke provided a formula relating the length $\lambda(I^{n+1}/JI^{n})$ to the difference $P_I(n)-H_I(n)$, where $I$ is an $\m$-primary ideal  in a $2$-dimensional Cohen-Macaulay local ring,  $J$ is a minimal reduction of $I$, \,$P_I(n)$ and $H_I(n)$ are respectively the  usual Hilbert-Samuel polynomial and the usual Hilbert-Samuel function of $I$ \cite{Huneke}. This formula was extended later by S. Huckaba to Cohen-Macaulay local rings of arbitrary dimension~$d$~\cite{Huckaba}. S.~Huckaba then established some formulas to compute the usual Hilbert coefficients of $I$,  and proved conditions in terms of  $e_1(I)$ for the associated graded ring   to be almost Cohen-Macaulay~\cite{Huckaba}.


If  $I$ is an $\m$-primary ideal in a Cohen-Macaulay local ring, the positivity of $e_1(I)$ can be observed from the well-known  Northcott's inequality
$$
e_1(I)\geq e_0(I)-\lambda(R/I)=\lambda(I/J),
$$
where $J$ is a minimal reduction of $I$.
By this inequality, one has that $e_1(I)=0$ if and only if $I$ is a complete intersection.
Furthermore, the ideal $I$ enjoys nice properties when equality holds in the above inequality.
Indeed, it was  shown that
$e_1(I)=\lambda(I/J)$ if and only if the reduction number of $I$ is at most 1, and in this case, the
associated graded ring $G$ is Cohen-Macaulay (see \cite{Huneke} and \cite{O}). 

This paper generalizes the above results to ideals that are not necessarily  $\m$-primary. 
In Section~2, we fix the notation and recall some basic concepts and facts that will be used throughout the paper. For an ideal  $I$  in a Noetherian local ring that has maximal analytic spread $\ell (I)=d={\rm dim}\,R$ and $G_d$ condition, we establish a formula to compute $e_1(\overline{I})$, where  $\overline{I}$ is an $1$-dimensional reduction of $I$ (see Section 2 for the definition of $\overline{I}$). We then give a condition in terms of $e_1(\overline{I})$ for
 the associated graded ring of $I$ to be almost Cohen-Macaulay.
 This result generalizes \cite[Theorem 3.1]{Huckaba}.  In Section~3, we  provide  a generalized
 version of \cite[Theorem 2.4]{Huckaba} relating the length
 $\lambda(I^{n+1}/JI^{n})$ to the difference $P_I(n)-H_I(n)$, where
 $I$ is an ideal in a $d$-dimensional Cohen-Macaulay local ring that
 satisfies $\ell(I)=d$, $G_d$ condition and $AN^-_{d-2}$,  $J$ is a
 general minimal reduction of $I$, \,$P_I(n)$ and $H_I(n)$ are respectively the
 generalized Hilbert-Samuel polynomial and the generalized Hilbert-Samuel function
 of $I$. As an application, we establish some formulas to compute the higher generalized
 Hilbert coefficients. In the last section, we apply our formula to prove a  generalized version of
 Northcott's inequality, and  recover the work of G. Colom$\acute{{\rm e}}$-Nin, C. Polini,
 B.  Ulrich and Y. Xie on the positivity of the generalized first  Hilbert coefficient $j_1(I)$.
 At the same time, we prove that, if equality holds in the generalized Northcott's inequality,
 the reduction number of $I$ is at most one and the  associated graded ring of $I$ is  Cohen-Macaulay,
 which generalizes the classical results of \cite{Huneke} and \cite{O}.

\section{Formulas for $e_1(\overline{I})$.}

\medskip

In  this paper,  we always assume that $(R, \m, k)$ is a Noetherian local ring of dimension $d$ with maximal ideal $\m$ and infinite residue field $k$ (we can enlarge the residue field to be infinite
 by replacing $R$ by $R(z)=R[z]_{\m R[z]}$, where $z$ is a variable over $R$).
Let $I$ be an $R$-ideal. We recall the concept of the generalized Hilbert-Samuel function of $I$.
Let $G={\rm gr}_{I} (R)=\oplus_{n=0}^{\infty}
I^n/I^{n+1}$ be the {\it associated graded ring}  of $I$.   As the homogeneous components of  $G$ may
not have  finite length,  one considers the $G$-submodule of
elements supported on $\m$:\,  $W:=\{\xi\in G\,|\, \exists\, t>0 \,\, {\rm such\, that}\,\, \xi \cdot \m^{t}=0\}
=H^0_{\m}(G)
=\oplus_{n=0}^{\infty} H^0_{\m}\,(I^n/I^{n+1})$.
Since $W$ is  a finite
graded module over ${\rm gr}_I(R) \otimes_R R/{\m}^{t}$ for some $t\geq 0$,
 its Hilbert-Samuel  function $H_W(n):=\sum_{i=0}^n \lambda(\Gamma_{\m}\,(I^i /I^{i +1}))$ is well defined.
 The {\it generalized Hilbert-Samuel function} of $I$ is defined to be: $H_{I}(n):=H_W(n)$ for every $n\geq 0$.

 The definition of generalized Hilbert-Samuel function was  introduced by H. Flenner, L. O' Carroll and W. Vogel in 1999  \cite[Definition 6.1.5]{FOV},
 and  studied later by C. Polini and Y. Xie \cite{PX} as well as G. Colom$\acute{{\rm e}}$ Nin, C. Polini, B.  Ulrich and Y. Xie \cite{NPUX}.
 Since ${\rm dim}_G\,W\leq {\rm dim}\,R=d$, \,  $H_{I}(n)$ is eventually a polynomial of degree  at most $d$
$$
P_{I}(n)=\sum_{i=0}^{d} (-1)^i j_i(I)
{{n+d-i}\choose{d-i}}.
$$
 C. Polini and Y. Xie \cite{PX} defined $P_{I}(n)$ to be the {\it generalized Hilbert-Samuel polynomial} of $I$ and $j_i(I),\, 0\leq i\leq d$,\, the {\it generalized Hilbert coefficients} of $I$. The normalized leading coefficient $j_0(I)$ is  called the {\it $j$-multiplicity} of $I$  (see  \cite{AM}, \cite{NU}, or \cite{PX}).

Recall that the Krull dimension of the special fiber ring
$G/{\m}G$ is called the
 {\it analytic spread} of $I$  and is denoted by $\ell(I)$. In general,
 ${\rm dim}_G\,W \leq \ell(I) \leq d$
 and equalities hold if and only if $\ell(I)=d$. Therefore $j_0(I)\neq 0$ if and only if $\ell(I)=d$ \cite{NU}.

If $I$ is $\m$-primary,  each homogeneous
component of $G$ has finite length, thus $W=G$ and the
generalized Hilbert-Samuel function coincides with the usual Hilbert-Samuel function; in particular, the generalized Hilbert coefficients $j_i(I)$, $0\leq i\leq d$, coincide with the usual Hilbert coefficients $e_i(I)$.

The definition of  generalized Hilbert coefficients is different from the one given by C. Ciuperc$\breve{{\rm a}}$ where he used the bigraded ring ${\rm gr}_{\m}(G)$ \cite{C}. Polini and Xie re-conciliated both approaches and proved that the generalized Hilbert coefficients  $j_0(I), \ldots, j_{d-2}(I)$ are preserved under a general hyperplane section \cite{PX}.

In this paper, we are going to use the tool of general elements 
to study the generalized Hilbert-Samuel function. We now recall this notion.
Let $I=(a_1, \ldots, a_{t})$ and write $x_i=\sum_{j=1}^{t} \lambda_{ij}a_j$ for $1\leq i\leq s$ and
$(\lambda_{ij})\in R^{s
t}$.  The elements  $x_1, \ldots, x_s$ form a {\it sequence of
general elements} in $I$ (equivalently  $x_1, \ldots, x_s$ are {\it
general} in $I$) if there exists a Zariski dense open subset $U$ of $k^{st}$
such that the image $(\overline{\lambda_{ij}})\in U$.
When $s=1$,  $x=x_1$ is said to be {\it general} in $I$.

Recall an ideal $J\subseteq I$ is called a {\it reduction}
 of $I$  if $JI^r=I^{r+1}$ for some non negative
 integer $r$.  The least such $r$ is denoted by $r_J(I)$.
 A reduction is {\it minimal} if it is
 minimal with respect to inclusion. The {\it reduction number} $r(I)$ of $I$ is defined as ${\rm min}\{r_J(I)\,|\,J \,\,{\rm is \,\,a \,\,mimimal\,\,reduction\,\,of \,\,}I\}$.
 Since $R$ has infinite residue field, the minimal number of
 generators $\mu(J)$ of any minimal reduction $J$  of $I$
 equals the analytic spread $\ell(I)$.   Furthermore, general $\ell(I)$ elements in $I$ form a minimal reduction $J$ whose
$r_J(I)$ coincides with the reduction number $r(I)$ (see \cite[2.2]{T2} or \cite[8.6.6]{SH}).
 One says that  $J$ is a {\it general minimal reduction} of $I$  if it is generated by $\ell(I)$ general elements in $I$.

 The ideal  $I$ is
  said to satisfy  {\it $G_{s+1}$ condition}  if for every
  $\p\in V(I)$ with ${\rm ht}\,\p=i\leq s$, the ideal
  $I_{\p}$ is generated by $i$ elements, i.e., $I_\p=(x_1,
  \ldots, x_{i})_\p$ for some $x_1, \ldots, x_{i}$ in $I$.

From now on,  we will assume $I$ has $\ell(I)=d$ and $G_d$ condition.
Let $J=(x_1, \ldots, x_{d})$, where $x_1, \ldots, x_{d}$ are general elements  in $I$, i.e., $J$ is  a general minimal reduction of $I$.
Set $J_{i}=(x_1, \ldots, x_{i})$, $0\leq i\leq d-1,$  \, $\overline{R}=R/J_{d-1}: I^{\infty}$, where $J_{d-1}: I^{\infty}=\{a\in R\,|\, \exists\, t>0 \,\, {\rm such\, that}\,\,  a \cdot I^{t}\subseteq J_{d-1}\}$, and    use
 $^{\tratto}$  to denote images in the quotient ring  $\overline{R}$. Then $\overline{R}$ is an $1$-dimensional Cohen-Macaulay local ring and $\overline{I}$ is $ \overline{\m}$-primary. Hence the generalized Hilbert-Samuel function $H_{\overline{I}}(n)$ and the generalized Hilbert-Samuel polynomial $P_{\overline{I}}(n)$ are the usual Hilbert-Samuel function and the usual Hilbert-Samuel polynomial of $\overline{I}$, respectively. Note $H_{\overline{R}}(I)$ and hence $P_{\overline{I}}(n)$ do not depend on choices of  general elements $x_1, \ldots, x_{d-1}$ in $I$ (see \cite{PX1}), and $P_{\overline{I}}(n)=e_0(\overline{I})(n+1)-e_1(\overline{I})$, where $e_0(\overline{I})=\lambda(\overline{R}/(\overline{x_d}))=j_0(I)$.
If $R$ is Cohen-Macaulay and $I$ is $\m$-primary, then $e_1(\overline{I})=e_1(I)$ (see for instance \cite[Proposition 1.2]{RV}).
But they are in general not the same.


We will show later in Theorem \ref{depthG} that $e_1(\overline{I})$ (like $e_1(I)$, see  \cite[Theorem 3.1]{Huckaba}) characterizes the depth of the associated graded ring $G$.
For ${\rm depth}(G)$, we mean the depth of the local ring $G_{M}$, where $M:=\m/I\oplus I/I^2 \oplus I^2/I^3\oplus \ldots$ denotes the maximal homogeneous ideal of $G$. Since ${\rm depth}(G)\leq {\rm dim}\,G={\rm dim}\,R=d$, $G$ is said to be {\it Cohen-Macaulay} if ${\rm depth}(G)= d$ and {\it almost Cohen-Macaulay} if ${\rm depth}(G)= d-1$.  The condition ${\rm depth}(G)\geq d-1$ is a useful one, especially when one considers questions about the behavior of $I^n$. It reduces greatly the computation of the generalized Hilbert coefficients (see Corollary 3.4 in Section 3).

Theorem \ref{depthG} is achieved from a formula  computing $e_1(\overline{I})$ (see Lemma \ref{Depth} in the following).  Since we do not have the  finite length  on $R/I$, to compare the length $\lambda(I^{n+1}/JI^n)$ with $\lambda(\overline{I}^{n+1}/\overline{J}\overline{I}^n)$, where $J$ is a general minimal reduction of $I$, we need the following lemma.

\begin{Lemma}\label{ABCD}
Let $D\subseteq B \subseteq A$ and $D\subseteq C \subseteq A$ be finite modules over $R$
such that $A/B$ and $C/D$ have finite lengths (while the lengths of $B/D$ and $A/C$ are not necessarily finite).
Then $$\lambda(A/B)+\lambda(B\cap C/D)=\lambda(C/D)+\lambda(A/B+C).$$
\end{Lemma}

\demo  By the exact sequences
$$
0{\rightarrow} B\cap C/D {\rightarrow} B/D \stackrel{\pi_1}{\rightarrow} B+C/C  {\rightarrow} 0,
$$
$$
0{\rightarrow} B+C/C \stackrel{i_1}{\rightarrow} A/C  {\rightarrow} A/B+C  {\rightarrow} 0,
$$
$$
0{\rightarrow} C/D  {\rightarrow} A/D \stackrel{\pi_2} {\rightarrow} A/C  {\rightarrow} 0,
$$
$$
0{\rightarrow} B/D \stackrel{i_2} {\rightarrow} A/D  {\rightarrow} A/B  {\rightarrow} 0,
$$
we have    the following long exact sequences
$$
0{\rightarrow} B\cap C/D {\rightarrow} H^0_{\m}(B/D)  {\rightarrow} H^0_{\m}(B+C/C) {\rightarrow} 0 {\rightarrow} H^1_{\m}(B/D) \stackrel{\widetilde{\pi_1}} {\rightarrow} H^1_{\m}(B+C/C) {\rightarrow} 0,
$$
$$
0{\rightarrow} H^0_{\m}(B+C/C) {\rightarrow} H^0_{\m}(A/C)  {\rightarrow} A/B+C \stackrel{\Delta_1} {\rightarrow} H^1_{\m}(B+C/C) \stackrel{\widetilde{i_1}}{\rightarrow} H^1_{\m}(A/C)   {\rightarrow} 0,
$$
$$
0{\rightarrow} C/D {\rightarrow} H^0_{\m}(A/D)  {\rightarrow} H^0_{\m}({A/C}) {\rightarrow} 0 {\rightarrow} H^1_{\m}(A/D) \stackrel{\widetilde{\pi_2}} {\rightarrow} H^1_{\m}(A/C) {\rightarrow} 0,
$$
$$
0 {\rightarrow} H^0_{\m}(B/D) {\rightarrow} H^0_{\m}(A/D)  {\rightarrow} A/B \stackrel{\Delta_2} {\rightarrow} H^1_{\m}(B/D) \stackrel{\widetilde{i_2}}{\rightarrow} H^1_{\m}(A/D)   {\rightarrow} 0,
$$
and the  commutative diagram
$$\begin{array}[c]{ccccccccc}
0 & {\rightarrow}& {\rm Im} (\Delta_2) & {\rightarrow}& H^1_{\m}(B/D) & \stackrel{\widetilde{i_2}}{\rightarrow} & H^1_{\m}(A/D) &{\rightarrow} &0\\
& & & & \downarrow \scriptstyle{ id} && \downarrow \scriptstyle{\widetilde{\pi_2}} &&\\
0 & {\rightarrow}& {\rm Ker} (\widetilde{i_1}\circ \widetilde{\pi_1}) & {\rightarrow}& H^1_{\m}(B/D) & \stackrel{\widetilde{i_1}\circ \widetilde{\pi_1}}{\rightarrow} & H^1_{\m}(A/C)&{\rightarrow} & 0
\end{array}
$$
with exact rows and isomorphic vertical maps $id$ and $\widetilde{\pi_2}$, hence ${\rm Im} (\Delta_2)\cong {\rm Ker} (\widetilde{i_1}\circ \widetilde{\pi_1})$. Since  ${\rm Ker} (\widetilde{i_1}\circ \widetilde{\pi_1})\cong {\rm Ker} (\widetilde{i_1})= {\rm Im} (\Delta_1)$, we have ${\rm Im} (\Delta_2)\cong {\rm Im} (\Delta_1)$. Now
$$
\lambda(A/B)+\lambda(B\cap C/D)=\lambda({\rm Im} (\Delta_2))+\lambda(H^0_{\m}(A/D))-\lambda(H^0_{\m}(B/D))+\lambda(B\cap C/D)
$$
$$
=\lambda({\rm Im} (\Delta_1))+\lambda(H^0_{\m}(A/D))-\lambda(H^0_{\m}(B/D))+\lambda(B\cap C/D)
$$
$$
=\lambda(A/B+C)+\lambda(H^0_{\m}(B+C/C)) - \lambda(H^0_{\m}(A/C)) +\lambda(H^0_{\m}(A/D))-\lambda(H^0_{\m}(B/D))+\lambda(B\cap C/D)
$$
$$
=\lambda(A/B+C)+\lambda(H^0_{\m}(B/D)) - \lambda(B\cap C/D) +\lambda(C/D)-\lambda(H^0_{\m}(B/D))+\lambda(B\cap C/D)
$$
$$
=\lambda(C/D)+\lambda(A/B+C).
$$
\QED

Applying Lemma \ref{ABCD}, we  obtain the following proposition.

\begin{Proposition}\label{Depth}
Let $I$ be an $R$-ideal with $\ell (I)=d$ and $G_d$ condition. For general elements $x_1, \ldots, x_d$ in $I$, set $J=(x_1, \ldots x_d)$, $J_{d-1}=(x_1, \ldots, x_{d-1})$, and $\overline{R}=R/J_{d-1}:I^{\infty}$ as above. Then for every $n\geq 0$, one has
\begin{itemize}
\item[(a)]
$\lambda(I^{n+1}/JI^{n})-\lambda[(J_{d-1}:I^{\infty})\cap I^{n+1}/(J_{d-1}:I^{\infty})\cap JI^{n}]=\Delta[ P_{\overline{I}}(n)-H_{\overline{I}}(n)].$
\smallskip
\item[(b)] $\sum_{n=0}^{\infty}[\lambda(I^{n+1}/JI^{n})-\lambda[(J_{d-1}:I^{\infty})\cap I^{n+1}/(J_{d-1}:I^{\infty})\cap JI^{n}]]=e_1(\overline{I}).$
\end{itemize}
\end{Proposition}

\demo  (a) For every $n\geq 0$, we have
$$\begin{array}[c]{ccc}
I^{n+1}+J_{d-1}:I^{\infty}& {\hookleftarrow} & I^{n+1}\\
\uparrow && \uparrow \\
JI^{n}+J_{d-1}:I^{\infty}& {\hookleftarrow} & JI^n
\end{array}$$
with $I^{n+1}+J_{d-1}:I^{\infty}/JI^{n}+J_{d-1}:I^{\infty}$ and $I^{n+1}/JI^{n}$ all having finite lengths.
By Lemma \ref{ABCD}, $$\lambda(I^{n+1}/JI^{n})=\lambda(I^{n+1}+J_{d-1}:I^{\infty}/JI^n+J_{d-1}:I^{\infty})+
\lambda([JI^n+J_{d-1}:I^{\infty}]\cap I^{n+1}/JI^n).$$ Since $[JI^n+J_{d-1}:I^{\infty}]\cap I^{n+1}/JI^n\cong (J_{d-1}:I^{\infty})\cap I^{n+1}/(J_{d-1}:I^{\infty})\cap JI^{n}$, we have
$$
\lambda(I^{n+1}/JI^{n})-\lambda[(J_{d-1}:I^{\infty})\cap I^{n+1}/(J_{d-1}:I^{\infty})\cap JI^{n}]=\lambda(\overline{I}^{n+1}/\overline{J}\overline{I}^{n})=\Delta[ P_{\overline{I}}(n)-H_{\overline{I}}(n)],
$$ where the latter equality follows from \cite[Theorem 2.4]{Huckaba}. Now (b) follows by (a) and \cite[Corollary 2.10]{Huckaba}.
\QED

\bigskip

We now recall some residual intersection properties.
Let $J_{i}=(x_1, \ldots, x_i)$, where $x_1,\ldots,x_i$ are elements in $I$. Define \, $J_{i}: I=
\{a\in R\,|\, a\cdot I \subseteq J_i\}.$
One says that  $J_i:I$ is an {\it  $i$-residual intersection } of~$I$
if $I_\p=(x_1, \ldots, x_{i})_\p$ for every $\p\in {\rm Spec}(R)$ with ${\rm dim}\,R_\p\leq i-1$.
An  $i$-residual intersection $J_i:I$ is called a {\it geometric $i$-residual intersection} of $I$
if, in addition,  $I_\p=(x_1, \ldots, x_{i})_\p$ for every
$\p\in V(I)$ with ${\rm dim}\,R_{\p}\leq i$.
It was shown that if $I$  satisfies   $G_s$ condition,
then for general elements $x_1, \ldots,  x_s$ in $I$ and each $0\leq i< s$,
the ideal $J_i : I$ is a geometric $i$-residual intersection of $I$,
and $J_s : I$ is an $s$-residual intersection of $I$ (see \cite{U} or \cite[Lemma~3.1]{PX}).

Assume  $R$ is   Cohen-Macaulay.
The ideal $I$ is  {\it $s$-weakly residually $(S_2)$} (respectively, has the {\it weak Artin-Nagata property} ${AN^-_s}$)
if for every $0\leq i\leq s$ and every geometric $i$-residual intersection $J_i:I$ of $I$
the quotient ring $R/J_i:I$ satisfies Serre's condition $(S_2)$ (respectively, is Cohen-Macaulay).

The notion of residual intersections was introduced by Artin and Nagata \cite{AN}
as a generalization of the concept of linkage to the case where the two ``linked" ideals
do not necessarily have the same height. The issue on the  Cohen-Macaulayness of residual intersections
 has been addressed in a series of results (for instance, \cite{Huneke1}, \cite{HVV}, \cite{HU} and \cite{U}),
 which require  either depth conditions  on all of the Koszul homology modules of $I$ such
 as  the ``strong Cohen-Macaulayness" or weaker ``sliding depth condition", or depth conditions on sufficiently many powers of $I$.

The following theorem generalizes  \cite[Theorem 3.1]{Huckaba} to ideals that are not necessarily
 $\m$-primary. Notice if $R$ is Cohen-Macaulay and $I$ is $\m$-primary,
  then $e_1(I)=e_1(\overline{I})$ (\cite[Proposition 1.2]{RV}), and  $I$ automatically satisfies $\ell(I)=d$,  $G_d$ condition, $(d-2)$-weakly residually $(S_2)$ as well as the weak Artin-Nagata property~$AN^-_{d-2}$.

\begin{Theorem} \label{depthG}
 Assume $R$ is Cohen-Macaulay. Let $I$ be an $R$-ideal which satisfies $\ell(I)=d$, $G_d$ condition and
 $(d-2)$-weakly residually $(S_2)$.
Then for   a general minimal reduction $J=(x_1, \ldots, x_d)$ of $I$,  the following two statements are equivalent:
\begin{itemize}
\item[(a)]
$\sum_{n=0}^{\infty}\lambda(I^{n+1}/JI^{n})=e_1(\overline{I}).$
\smallskip
\item[(b)] For every $n\geq 0$, \,$J_{d-1}\cap I^{n+1}=J_{d-1}I^n$,
where $J_{d-1}=(x_1, \ldots, x_{d-1})$ defined as before.
\end{itemize}
Furthermore, if $I$ satisfies $AN^-_{d-2}$,  then  (a) or (b) is equivalent to that ${\rm depth}
(G)\geq d-1$. 
\end{Theorem}

\demo First if $d=1$, then $\overline{R}=R/0:I^{\infty}=R/0:I$ is an 1-dimensional Cohen-Macaulay local
ring and $\overline{I}$ is $\overline{\m}$-primary. By $(0:I)\cap I=0$ and \cite[Theorem 3.1]{Huckaba}, one has
$$
e_1(\overline{I})=\sum_{n=0}^{\infty}\lambda(\overline{I}^{n+1}/\overline{J}
\overline{I}^{n})=\sum_{n=0}^{\infty}\lambda(I^{n+1}/JI^{n}).
$$

Assume $d\geq 2$.  Since $I$ satisfies $\ell(I)=d$ and $G_d$ condition,  one has that $J_{i}: I$
is a geometrically $i$-residual intersection of $I$, where
 $J_{i}=(x_1, \ldots, x_{i})$, $0\leq i \leq d-1$ \cite{PX}.
Furthermore, since $I$ is $(d-2)$-weakly residually $(S_2)$, for each $0\leq i \leq d-1$, one has that
$R/J_i:I$ has no embedded associated prime ideals
and thus $J_i:I=J_i: x_{i+1}$. Note that ${\rm depth}(R/J_i:x_{i+1})={\rm depth}(R/J_i:I)\geq 2$.
We will show
${\rm depth}(R/J_i)\geq 1$ by induction on $i$.
The case  $i=0$ (i.e., $J_0=(0)$) is clear.  Assume $1\leq i\leq d-1$ and ${\rm depth}(R/J_{i-1})\geq 1$, then
by the exact sequence
$$
0\rightarrow R/J_{i-1}:x_{i}\rightarrow R/J_{i-1}\rightarrow R/J_i \rightarrow 0,
$$
one has that ${\rm depth}(R/J_i)\geq {\rm Min}\{{\rm depth}(R/J_{i-1}:x_i) -1, \, {\rm depth}(R/J_{i-1})\}\geq 1$.
We claim that
$
(J_{d-1}: I)\cap I =J_{d-1}.
$
Indeed, since $J_{d-1}\subseteq (J_{d-1}: I)\cap I$, we just need to show that $(x_1, \ldots, x_{d-1})_{\p}= (J_{d-1}: I)_{\p}\cap I_{\p}$ for every $\p\in {\rm Ass}(R/J_{d-1})$, which  follows by the fact that for every $\p\in {\rm Ass}(R/J_{d-1})$,
since ${\rm height}\, \p\leq d-1$, then either $I_{\p}=R_{\p}$ or $I_{\p}=(J_{d-1})_{\p}$.

Now for $n\geq 0$, $(J_{d-1}: I)\cap I^{n+1}=J_{d-1}\cap I^{n+1}$ and $(J_{d-1}: I)\cap JI^{n}=J_{d-1}\cap JI^{n}$.
Therefore if  (b) is true, then for $n\geq 0$,
$$
\lambda[(J_{d-1}: I)\cap I^{n+1}/(J_{d-1}: I)\cap JI^{n}]=\lambda[J_{d-1}\cap I^{n+1}/J_{d-1}\cap JI^{n}]
$$
$$
=\lambda[J_{d-1}I^n/J_{d-1}\cap JI^{n}]=0.
$$
And  (a) follows by Proposition \ref{Depth} (b) and the fact that $R/J_{d-1}: I$ is Cohen-Macaulay and thus $J_{d-1}: I^{\infty}=J_{d-1}: I$.

Assume  (a). By Proposition \ref{Depth} (b), for every $n\geq 0$,\, $\lambda[(J_{d-1}: I)\cap I^{n+1}/(J_{d-1}: I)\cap JI^{n}]=0$. Hence
$$
J_{d-1}\cap I^{n+1}=(J_{d-1}: I)\cap I^{n+1}=(J_{d-1}: I)\cap JI^{n}=J_{d-1}\cap JI^{n}.
$$
We use induction on $n$ to prove that for every $n\geq 0,$ $J_{d-1}\cap I^{n+1}=J_{d-1}I^n$.
This is clear if $n=0$. Assume $n\geq 1$ and $J_{d-1}\cap I^{n}=J_{d-1}I^{n-1}$. Then (b) follows by the following equalities:
\begin{eqnarray*}
J_{d-1}\cap I^{n+1}&=&J_{d-1}\cap JI^{n} \\
&=&J_{d-1}\cap (J_{d-1}I^{n}+x_{d}I^{n})\\
&=& J_{d-1}I^{n}+J_{d-1}\cap x_{d}I^{n}\\
&=& J_{d-1}I^{n}+x_{d}[(J_{d-1}: x_{d})\cap I^{n}]\\
&=&J_{d-1}I^{n}+x_{d}[J_{d-1}\cap I^{n}] \\ 
&= &J_{d-1}I^{n}+x_{d}J_{d-1} I^{n-1} \\
&= &J_{d-1}I^{n}.
\end{eqnarray*}

Finally assume $I$ satisfies $AN^-_{d-2}$, we will show that  (b) is equivalent to that ${\rm depth} (G)\geq d-1$.
Set  $\delta(I)=d-g$, where ${\rm ht}\,I=g$. We use the induction on $\delta$.  If  $\delta=0$, the assertion follows because (b) is equivalent to that $x_1^{*}, \ldots, x_{d-1}^{*}$ form a regular sequence on  $G$
(see \cite[Proposition~2.6]{VV}), and the latter is equivalent to that ${\rm depth} (G)\geq d-1$.
Thus we may assume  $\delta(I)\geq 1$ and the theorem holds for smaller values of $\delta(I)$. In particular, $d\geq g+1$. Since $x_1^{*}, \ldots, x_g^{*}$ form a  regular sequence on  $G$, we may factor out $x_1, \ldots, x_g$ to assume $g=0$. Now $d=\delta(I) \geq 1$.
Set   $S=R/0: I$. Then $S$ is Cohen-Macaulay  since $I$ satisfies $AN^-_{d-2}$.
Note ${\rm dim}\, S={\rm dim}\,R=d$,  ${\rm grade}\,(IS)\geq 1$, $IS$ still satisfies  $G_{d}$ condition, $AN^-_{d-2}$,  $\ell(IS)=\ell(I)=d$ (see for instance \cite{PX}).
 Since $I\cap (0: I)=0$, there is a graded exact sequence
\begin{equation}\label{eq11}
0\rightarrow 0: I \rightarrow G \rightarrow {\rm gr}_{IS}(S)\rightarrow 0.
\end{equation}
Since ${\rm depth}(0: I)\geq d$,  one has that ${\rm depth} (G)\geq  d-1 \Leftrightarrow {\rm depth}({\rm gr}_{IS}(S))\geq d-1$.
 We claim  that (b) is equivalent to  $J_{d-1}S\cap I^{n+1}S=J_{d-1}I^n S$ for every $n\geq 0$.
Indeed, if (b) holds, then clearly $J_{d-1}S\cap I^{n+1}S=J_{d-1}I^n S$ for every $n\geq 0$. On the other hand,
if $J_{d-1}S\cap I^{n+1}S=J_{d-1}I^n S$ for every $n\geq 0$, then
  $$
  J_{d-1}\cap I^{n+1}\subseteq J_{d-1}I^n+(0: I)\cap I^{n+1}=J_{d-1}I^n,
  $$
again by   $I\cap 0: I=0$. 
We are done by induction hypothesis since  $\delta(IS)=d-{\rm grade}\,(IS)<d=\delta(I)$.
 
 \QED

\section{Formulas for $j_i(I)$, $1\leq i\leq d$. }

 In this section we will provide a formula relating the length $\lambda(I^{n+1}/JI^{n})$ to the difference $P_I(n)-H_I(n)$, where $I$ is an ideal with $\ell(I)=d$, $G_d$ condition and  $AN^-_{d-2}$,
 $J$ is a general minimal reduction of $I$, $P_I(n)$ and $H_I(n)$ are  the
 generalized Hilbert-Samuel polynomial and the generalized Hilbert-Samuel function of $I$.
This formula generalizes \cite[Theorem 2.4]{Huckaba}

\begin{Theorem} \label{Theorem2}
Assume $R$ is Cohen-Macaulay.  Let  $I$ be  an $R$-ideal with  $\ell(I)=d$,  $G_d$ condition and  $AN^-_{d-2}$.
Then for  a general minimal reduction $J=(x_1, \ldots, x_d)$ of $I$, one has that
  for all $n\geq 0$,
$$
\lambda(I^{n+1}/JI^{n})+\omega_{n}(J,I)=\Delta^d[ P_{I}(n)-H_{I}(n)],
$$
where $\omega_0(J, I)=\lambda(R/J_{d-1}:I+I)-\lambda[H^0_{\m}(R/I)]$,
and for $n\geq 1$,
\begin{eqnarray*}
\omega_{n}(J,I)&=& \Delta^{d-1}[\lambda({\widetilde{K}^0_{n-1}})]+\Delta^{d-2}[\lambda({\widetilde{K}^1_{n-1}})]
+\ldots + \Delta[\lambda(\widetilde{K}^{d-2}_{n-1})]\\
&+&\Delta^{d-2}[\lambda(\widetilde{L}^0_{n})-\lambda(L^0_{n})+\lambda(N^0_n)]\\
&+& \Delta^{d-3}[\lambda(\widetilde{L}^1_{n})-\lambda(L^1_{n})+\lambda(N^1_n)]+\ldots\\
&+& \Delta^0[\lambda(\widetilde{L}^{d-2}_{n})-\lambda(L^{d-2}_{n})+\lambda(N^{d-2}_n)]\\
&-& \lambda[(J_1:I)\cap I^{n+1}/(J_1:I)\cap JI^n]-\lambda[(J_2:I)\cap I^{n+1}+(J_1:I)/(J_2:I)
\cap JI^n+(J_1:I)]\\
&& -\ldots - \lambda[(J_{d-1}:I)\cap I^{n+1}+(J_{d-2}:I)/(J_{d-1}:I)\cap JI^n+(J_{d-2}:I)]\\
&&-(-1)^n\binom{d-1}{n}\beta,
\end{eqnarray*}
$\binom{d-1}{n}:=0$ if $n\geq d$, and for $0\leq i\leq d-2$,
$$
\widetilde{K}^i_{n-1}=I^{n+1}:x_1/J_i:I+I^n,
$$
$$
\widetilde{L}^i_{n}=J_{i+1}\cap I^n/[J_{i}\cap I^{n}+J_{i+1}\cap I^{n+1}+x_{i+1} I^{n-1}],
$$
$$L^i_n=((J_i:I)\cap I^n+I^{n+1}):_{(J_{i+1}:I)\cap I^{n}} \m^{\infty}/[(J_i:I)\cap I^n+(J_{i+1}:I)\cap I^{n+1}+x_{i+1}(((J_i:I)\cap I^{n-1}+I^{n}):_{I^{n-1}}\m^{\infty})],$$
 $$N^i_n=((J_{i+1}:I)\cap I^n+I^{n+1}):_{I^n}\m^{\infty}/[(J_{i+1}:I)\cap I^n+((J_i:I)\cap I^n+I^{n+1}):_{I^n}\m^{\infty}],$$
 $$
 \beta=\lambda(H^0_{\m}(R/I))-\lambda(H^0_{\m}(R/0:I+I)).
 $$
\end{Theorem}

\demo
Recall for each $0\leq i\leq d-1$, \, $J_i:I$ is a geometric $i$-residual intersection of $I$, where $J_i=(x_1, \ldots, x_i)$. Set
$R^i=R/J_i:I$ and $G^i={\rm gr}_{IR^i}(R^i)$.
 Then $[G^{0}]_{0}=R/(0:I+I)$ and
 $[G^{0}]_{n}=[G]_{n}$ for every $n\geq 1$. Hence
\begin{eqnarray*}
&& \Delta[H_{I}(0)]=\lambda(H^0_{\m}(R/I))\\
 &=&\lambda(H^0_{\m}(R/0:I+I))+[\lambda(H^0_{\m}(R/I))-\lambda(H^0_{\m}(R/0:I+I))]\\
 &=&\Delta[H_{IR^0}(0)]
 +\beta,
\end{eqnarray*}
with $\beta$  defined above,  and  $\Delta[H_{I}(n)]=\Delta[H_{IR^{0}}(n)]$ for $n\geq 1$.
 Therefore we have that for $n\geq 0$,
  \begin{equation}\label{eq1}
 \Delta^d[H_I(n)]=\Delta^d[H_{IR^0}(n)]+(-1)^n \binom{d-1}{n}\beta, 
\end{equation}
with the  binomial coefficient  $\binom{d-1}{n}=0$ if $n>d-1$.

 We use induction on $d$ to prove the theorem. First assume $d=1$.
If $n=0$, one has
\begin{eqnarray*}
&&\lambda(I/J)+\omega_0(J, I)\\
&=&\lambda(IR^{0}/JR^{0})+\lambda(R/0:I+I)-\lambda[H^0_{\m}(R/I)]\\
&=&\Delta[P_{IR^{0}}(0)-H_{IR^{0}}(0)]+\lambda(R/0:I+I)-\lambda[H^0_{\m}(R/I)]\\
&=&\Delta[P_{IR^{0}}(0)]-\lambda(R/0:I+I)+\lambda(R/0:I+I)-\lambda[H^0_{\m}(R/I)]\\
&=&\Delta[P_{I}(0)-H_{I}(0)],
\end{eqnarray*}
where the second equality follows from \cite[Theorem 2.4]{Huckaba} since
$R^{0}$ is an 1-dimensional Cohen-Macaulay local ring and $IR^{0}$ is $\m R^0$-primary, and the third equality follows from $\Delta[P_{IR^{0}}(n)]=\Delta[P_{I}(n)]=j_0(I)$ for every $n\geq 0$.

If $n\geq 1$, since $\omega_n(J, I)=0$, one  has
$$
\lambda(I^{n+1}/JI^n)+\omega_n(J, I)=\lambda(I^{n+1}R^{0}/JI^nR^{0})
$$
$$
=\Delta[P_{IR^{0}}(n)-H_{IR^{0}}(n)]
=\Delta[P_{I}(n)-H_{I}(n)].
$$

Now assume $d\geq 2$ and the assertion holds for $d-1$.
By the proof of Proposition \ref{Depth},
$$\lambda(I^{n+1}/JI^{n})-\lambda[(J_{d-1}:I)\cap I^{n+1}/(J_{d-1}:I)\cap JI^n]=\lambda(I^{n+1}R^{d-1}/JI^{n}R^{d-1})
$$
$$
=
\Delta[ P_{IR^{d-1}}(n)-H_{IR^{d-1}}(n)]=\Delta^d[ P_{I}(n)]-\Delta[H_{IR^{d-1}}(n)], 
$$
by the fact that $\Delta[ P_{IR^{d-1}}(n)]=\Delta^d[ P_{I}(n)]=j_0(I)$. 
If $n=0$, one has 
$$\lambda[(J_{d-1}:I)\cap I/(J_{d-1}:I)\cap J]=\lambda(J_{d-1}/J_{d-1})=0,$$
and therefore
\begin{eqnarray*}
&&\lambda(I/J)+\omega_{0}(J,I)\\
&=&\Delta^d[ P_{I}(0)]-\Delta[H_{IR^{d-1}}(0)]+\lambda(R/J_{d-1}:I+I)-\lambda[H^0_{\m}(R/I)]\\
&=&\Delta^d[ P_{I}(0)]-\lambda(R/J_{d-1}:I+I)+\lambda(R/J_{d-1}:I+I)-\lambda[H^0_{\m}(R/I)]\\
&=&\Delta^d[ P_{I}(0)]-\lambda[H^0_{\m}(R/I)]\\
&=&\Delta^d[ P_{I}(0)-H_I(0)].
\end{eqnarray*}

Let $n\geq 1$. We have the following exact sequences for every $n\geq 1$:
 $$
 0\rightarrow K^0_{n-1} \rightarrow H^0_{\m}([G^0]_{n-1}) \stackrel{x_1^*} \rightarrow H^0_{\m}([G^0]_{n})\rightarrow H^0_{\m}([G^0]_{n})/x_1^*H^0_{\m}([G^0]_{n-1}) \rightarrow 0,
 $$
 $$
0 \rightarrow L^0_{n} \rightarrow H^0_{\m}([G^0]_{n})/x_1^*H^0_{\m}([G^0]_{n-1})  \rightarrow H^0_{\m}
([G^{1}]_n) \rightarrow
N^0_n
 \rightarrow 0,
 $$
 where 
 $$K^0_{n-1}=[((0:I)\cap I^n+I^{n+1}):_{I^{n-1}}x_1]\cap [((0:I)\cap I^{n-1}+I^{n}):_{I^{n-1}}\m^{\infty}]/((0:I)\cap I^{n-1}+I^n),$$
 $$L^0_n=((0:I)\cap I^n+I^{n+1}):_{(J_1:I)\cap I^{n}} \m^{\infty}/[(0:I)\cap I^n+(J_1:I)\cap I^{n+1}+x_1(((0:I)\cap I^{n-1}+I^{n}):_{I^{n-1}}\m^{\infty})],$$
 $$N^0_n=((J_1:I)\cap I^n+I^{n+1}):_{I^n}\m^{\infty}/[(J_1:I)\cap I^n+((0:I)\cap I^n+I^{n+1}):_{I^n}\m^{\infty}].$$
 Note  $((0:I)\cap I^n+I^{n+1}):_{I^{n-1}}x_1/((0:I)\cap I^{n-1}+I^n)$  has finite
 length because $G$ is Cohen-Macaulay on the punctured spectrum by \cite[Theorem 3.1]{JU}.
Hence  $$K^0_{n-1}=((0:I)\cap I^n+I^{n+1}):_{I^{n-1}}x_1/((0:I)\cap I^{n-1}+I^n).$$ 
 Therefore
 \begin{eqnarray*}
 \Delta^d[H_{IR^0}(n)]&=&\Delta^{d-2}[\lambda[H^0_{\m}([G^0]_{n})]-\lambda[H^0_{\m}([G^0]_{n-1})]] \\
 &=&\Delta^{d-2}[\lambda[H^0_{\m}([G^0]_{n})/x_1^*H^0_{\m}([G^0]_{n-1})]-\lambda(K^0_{n-1})]\\
 &=& \Delta^{d-2}[\lambda[H^0_{\m}([G^{1}]_n)]]+\Delta^{d-2}
 [\lambda(L^0_{n})]-\Delta^{d-2}
 [\lambda(N^0_{n})]
 -\Delta^{d-2}
 [\lambda(K^0_{n-1})]\\
 &=&\Delta^{d-1}[H_{IR^{1}}(n)]
 +\Delta^{d-2}
 [\lambda(L^0_{n})-
 \lambda(N^0_{n})
 -\lambda(K^0_{n-1})].
\end{eqnarray*}
By  Lemma  \ref{ABCD},  the induction hypothesis,  and the above equality,
\begin{eqnarray*}
&&\lambda(I^{n+1}/JI^{n})=\lambda(I^{n+1}R^0/JI^{n}R^0)\\
&=&\lambda(I^{n+1}R^{1}/JI^{n}R^{1})+\lambda[(J_{1}:I)\cap I^{n+1}/(J_{1}:I)\cap JI^n]\\
&=&
\Delta^{d-1}[ P_{IR^{1}}(n)-H_{IR^{1}}(n)]-\omega_n(JR^{1}, IR^{1})+\lambda[(J_{1}:I)\cap I^{n+1}/(J_{1}:I)\cap JI^n]\\
&=&\Delta^d[P_{IR^0}(n)]-\Delta^d[H_{IR^0}(n)]+\Delta^{d-2}
 [\lambda(L^0_{n})-
 \lambda(N^0_{n})
 -\lambda(K^0_{n-1})]\\
 && - \omega_n(JR^{1}, IR^{1})+\lambda[(J_{1}:I)\cap I^{n+1}/(J_{1}:I)\cap JI^n]\\
&=&\Delta^d[P_{IR^0}(n)-H_{IR^0}(n)]\\
&& -[\omega_n(JR^{1}, IR^{1})
 + \Delta^{d-2}
 [\lambda(K^0_{n-1})]+ \Delta^{d-2}[-\lambda(L^0_{n})+
 \lambda(N^0_{n})]-\lambda[(J_{1}:I)\cap I^{n+1}/(J_{1}:I)\cap JI^n]].
\end{eqnarray*}

We claim  that for every $n\geq 1$,
$$
\lambda(K^0_{n-1})
=\Delta[\lambda(\widetilde{K}^0_{n-1})]+\lambda(\widetilde{L}^0_{n}),
$$
where
$$
\lambda(\widetilde{K}^0_{n-1})=\lambda(I^{n+1}:x_1/0:I+I^n),
$$
$$
\lambda(\widetilde{L}^0_{n})=\lambda[(x_1)\cap I^n/(x_1)\cap I^{n+1}+x_1 I^{n-1}],
$$
since $(0:I)\cap I^n=0$ for $n\geq 1$. 
This follows by the following equalities: 
\begin{eqnarray*}
&& \lambda(K^0_{n-1})=\lambda[I^{n+1}:_{I^{n-1}}x_1/(0:I)\cap I^{n-1}+I^n]\\
&=&\lambda[(I^{n+1}:_{I^{n-1}}x_1)+0:I/I^n+0:I]\\
&=&\lambda[I^{n+1}:x_1/I^n+0:I]-\lambda[I^{n+1}:x_1/(I^{n+1}:_{I^{n-1}}x_1)+0:I]\\
&=&\Delta[\lambda(\widetilde{K}^0_{n-1})]+\lambda[I^{n}:x_1/0:I+I^{n-1}]-\lambda[(I^{n+1}:x_1)+I^{n-1}/0:I+I^{n-1}]\\
&=&\Delta[\lambda(\widetilde{K}^0_{n-1})]+\lambda[I^{n}:x_1/(I^{n+1}:x_1)+I^{n-1}]\\
&=&\Delta[\lambda(\widetilde{K}^0_{n-1})]+\lambda[(x_1)\cap I^n/(x_1)\cap I^{n+1}+x_1 I^{n-1}]\\
&=&\Delta[\lambda(\widetilde{K}^0_{n-1})]+\lambda(\widetilde{L}^0_{n}).
\end{eqnarray*}
Now
\begin{eqnarray*}
&&\omega_n(JR^1,IR^1)+\Delta^{d-2}
 [\lambda(K^0_{n-1})]+\Delta^{d-2}
 [-\lambda(L^0_{n})+\lambda(N^0_{n})]-\lambda[(J_{1}:I)\cap I^{n+1}/(J_{1}:I)\cap JI^n]\\
&=&\omega_n(JR^1,IR^1)+\Delta^{d-1}
 [\lambda(\widetilde{K}^0_{n-1})]+\Delta^{d-2}
 [\lambda(\widetilde{L}^0_{n})-\lambda(L^0_{n})+\lambda(N^0_{n})]-\lambda[(J_{1}:I)\cap I^{n+1}/(J_{1}:I)\cap JI^n]\\
&=&\omega_n(JR^0,IR^0).
\end{eqnarray*}
Therefore by equation (\ref{eq1}),  we have 
$$
\lambda(I^{n+1}/JI^{n})=
=\Delta^d[P_{IR^0}(n)-H_{IR^0}(n)]-\omega_n(JR^0,IR^0).
$$
$$
=\Delta^d[P_{I}(n)-H_{I}(n)]-[\omega_n(JR^0,IR^0)-(-1)^n\binom{d-1}{n}\beta]=\Delta^d[P_{I}(n)-H_{I}(n)]-\omega_n(J,I).
$$
\QED

\bigskip

The following Lemma is inspired by \cite[Proposition 2.9]{Huckaba}.

\begin{Lemma}
Let $I$ be  an $R$-ideal. Then
$$
\sum_{n={i-1}}^{\infty}\binom{n}{i-1}\Delta^d[P_I(n)-H_I(n)]=j_i(I) \,\, {\rm for}\,\, 1\leq i \leq d.
$$
\end{Lemma}

\bigskip

By Theorem \ref{Theorem2} and the above lemma, we obtain  formulas to compute the generalized Hilbert coefficients.

\begin{Corollary}\label{j1}
Assume  $R$ is Cohen-Macaulay. Let  $I$ be  an $R$-ideal with  $\ell(I)=d$,  $G_d$
condition and  $AN^-_{d-2}$. Then for a general minimal reduction  $J=(x_1, \ldots, x_d)$  of $I$,
 one has
 $$
 \sum_{n={i-1}}^{\infty}\binom{n}{i-1}[\lambda(I^{n+1}/JI^n)+\omega_n(J,I)]=j_i(I) \,\, {\rm for}\,\, 1\leq i \leq d.
 $$
 In particular, if $d=1$,
 $$
 j_1(I)
=\sum_{n={0}}^{\infty}\lambda(I^{n+1}/JI^{n})+\lambda(R/0:I+I)-\lambda[H^0_{\m}(R/I)],
$$
and if $d\geq 2$,
 \begin{eqnarray*}
j_1(I)&
=&\sum_{n={0}}^{\infty}\lambda(I^{n+1}/JI^{n})+\lambda(R/J_{d-1}:I+I)-\lambda[H^0_{\m}(R/0:I+I)] \\
&+& \Delta^{d-2}[\lambda(L_0^0)-\lambda(N_0^0)]+\ldots +\Delta[\lambda(L_0^{d-3})-\lambda(N_0^{d-3})]\\
&+&\sum_{n=1}^{\infty}[\lambda(\widetilde{L}^{d-2}_{n})-\lambda(L^{d-2}_{n})+\lambda(N^{d-2}_n)]\\
&-&\sum_{n=1}^{\infty}\lambda((J_1:I)\cap I^{n+1}/(J_1:I)\cap JI^n)-\ldots \\
 &-& \sum_{n=1}^{\infty}\lambda[(J_{d-1}:I)\cap I^{n+1}+J_{d-2}:I/(J_{d-1}:I)\cap JI^n+J_{d-2}:I].
\end{eqnarray*}
\end{Corollary}

\demo
If $d=1$, by  Theorem \ref{Theorem2}, one has $\omega_0(J, I)=\lambda(R/0:I+I)-\lambda[H^0_{\m}(R/I)]$ and 
$\omega_n(J, I)=0$ for $n\geq 1$. Hence
 $$
 j_1(I)=\sum_{n={0}}^{\infty}[\lambda(I^{n+1}/JI^n)+\omega_n(J,I)]
=\sum_{n={0}}^{\infty}\lambda(I^{n+1}/JI^{n})+\lambda(R/0:I+I)-\lambda[H^0_{\m}(R/I)]. 
$$
Assume $d\geq 2$. Then 
\begin{eqnarray*}
j_1(I)&=&\sum_{n={0}}^{\infty}[\lambda(I^{n+1}/JI^n)+\omega_n(J,I)]\\
&=&\sum_{n={0}}^{\infty}\lambda(I^{n+1}/JI^{n})+\lambda(R/J_{d-1}:I+I)-\lambda[H^0_{\m}(R/I)] \\
&+& \Delta^{d-2}[\lambda(L_0^0)-\lambda(N_0^0)]+\ldots +\Delta[\lambda(L_0^{d-3})-\lambda(N_0^{d-3})]\\
&+&\sum_{n=1}^{\infty}[\lambda(\widetilde{L}^{d-2}_{n})-\lambda(L^{d-2}_{n})+\lambda(N^{d-2}_n)]\\
&-&\sum_{n=1}^{\infty}\lambda((J_1:I)\cap I^{n+1}/(J_1:I)\cap JI^n)-\ldots \\
 &-& \sum_{n=1}^{\infty}\lambda[(J_{d-1}:I)\cap I^{n+1}+J_{d-2}:I/(J_{d-1}:I)\cap JI^n+J_{d-2}:I]\\
 &-& \beta \big[\sum_{n=0}^{d-1}(-1)^n\binom{d-1}{n}\big] + \beta,
\end{eqnarray*}
which is equal to the desired result since
$\sum_{n=0}^{d-1}(-1)^n\binom{d-1}{n}=0$  and $$\beta=\lambda(H^0_{\m}(R/I))-\lambda(H^0_{\m}(R/0:I+I)).$$
\QED

\smallskip

\begin{Corollary}
Assume $R$ is  Cohen-Macaulay. Let  $I$ be an $R$-ideal with  $\ell(I)=d$,  $G_d$ condition and  $AN^-_{d-2}$. If
${\rm depth}(G)\geq d-1$ and ${\rm depth}(G/H^0_{\m}(G))=d$, then  for a general minimal reduction  $J=(x_1, \ldots, x_d)$  of $I$, one has
$$
j_1(I)
=\sum_{n={0}}^{\infty}\lambda(I^{n+1}/JI^{n})+\lambda(R/J_{d-1}:I+I)-\lambda[H^0_{\m}(R/H+I)],
$$
where $H=0$ if $d=1$, or $H=0:I$ if  $d\geq 2$.
\end{Corollary}

\section{Generalized Northcott's inequality}

As an application of Corollary \ref{j1}, we obtain the following generalized Northcott's inequality.

  \begin{Theorem}\label{North}
  Assume $R$ is  Cohen-Macaulay. Let  $I$ be  an $R$-ideal with $\ell(I)=d$,  $G_d$ condition
  and weakly $(d-2)$-residually $(S_2)$.  Then for a general minimal reduction  $J=(x_1, \ldots, x_d)$ of $I$, one has the following generalized Northcott's inequality:
   $$
   j_1(I)\geq \lambda(I/J)+\lambda[R/J_{d-1}:I+(J_{d-2}:I+I):\m^{\infty}].
   $$
   \end{Theorem}
   
\demo Set $S=R/J_{d-2}:I$, where $J_{d-2}=(x_1, \ldots, x_{d-2})$.
 Then $j_1(I)=j_1(IS)$, \,$IS$ satisfies  $\ell(IS)=2={\rm dim}\,S$,  $G_2$ condition and $AN^-_0$ (see \cite{PX1} and \cite{PX}).
  By Corollary \ref{j1}, we have
  \begin{eqnarray*}
 & & j_1(I)=j_1(IS)\\
&=&\sum_{n=0}^{\infty}\lambda(I^{n+1}S/JI^{n}S)-\sum_{n=1}^{\infty}
\lambda[(x_{d-1}S:IS)\cap I^{n+1}S/(x_{d-1}S:IS)\cap JI^{n}S]\\
 & + &\lambda(S/x_{d-1}S: IS+IS)-\lambda(H^0_{\m}(S/IS))\\
& + &\sum_{n=1}^{\infty}\big[\lambda[(x_{d-1}S)\cap I^nS/(x_{d-1}S)\cap I^{n+1}S+x_{d-1} I^{n-1}S]\\
&-&\lambda[I^{n+1}S:_{(x_{d-1}S:IS)\cap I^{n}S} \m^{\infty}/(x_{d-1}S:IS)\cap I^{n+1}S+x_{d-1}(I^{n}S:_{I^{n-1}S}\m^{\infty})]\big]\\
&+&
\sum_{n=1}^{\infty} \lambda\big[((x_{d-1}S:IS)\cap I^nS+I^{n+1}S):_{I^nS}\m^{\infty}/(x_{d-1}S:IS)\cap I^nS+I^{n+1}S:_{I^nS}\m^{\infty}\big]\\
& \geq & \lambda(I/J)+\lambda[R/J_{d-1}:I+(J_{d-2}:I +I):\m^{\infty}].
\end{eqnarray*}
This follows by the following inequalities. First
 \begin{eqnarray*}
&& \sum_{n=0}^{\infty}\lambda(I^{n+1}S/JI^{n}S)-\sum_{n=1}^{\infty}
\lambda[(x_{d-1}S:IS)\cap I^{n+1}S/(x_{d-1}S:IS)\cap JI^{n}S]\\
& =& \lambda(IS/JS)+ \sum_{n=1}^{\infty}\big[\lambda(I^{n+1}S/JI^{n}S)-
\lambda[(x_{d-1}S:IS)\cap I^{n+1}S/(x_{d-1}S:IS)\cap JI^{n}S]\big]\\
&=& \lambda(I/J) + \sum_{n=1}^{\infty}\lambda[I^{n+1}S/JI^{n}S+ (x_{d-1}S:IS)\cap I^{n+1}S]\\
&\geq & \lambda(I/J),
 \end{eqnarray*}
where the second equality follows by Lemma \ref{ABCD} and 
$\lambda(IS/JS)=\lambda(I/(J_{d-1}:I)\cap I +J)=\lambda(I/J)$. 

Next, because ${\rm depth}(S/x_{d-1}S)\geq 1$ (see \cite{U}),  for  every $\p\in {\rm Ass}(S/x_{d-1}S)$, one has that $\p$ is not maximal and $IS_{\p}=(x_{d-1})S_{\p}$. Hence  $(x_{d-1}S: IS)\cap (JS: \m^{\infty})_{\p}=x_{d-1}S_{\p}$ for every $\p\in {\rm Ass}(S/x_{d-1}S)$, which yields that
$(x_{d-1}S: IS)\cap (JS: \m^{\infty})=x_{d-1}S$.
Therefore 
 $$(x_{d-1}S: IS+JS)\cap (JS: \m^{\infty})= JS + (x_{d-1}S: IS)\cap (JS: \m^{\infty})=JS. $$ 
 Since $\lambda(I/J)<\infty$ and $(x_{d-1}S: IS)\cap IS=(x_{d-1})S$, by Lemma \ref{ABCD}, one has 
 \begin{eqnarray*}
&& \lambda(S/x_{d-1}S: IS+IS)-\lambda(H^0_{\m}(S/IS))\\
&=& \lambda(S/x_{d-1}S: IS+JS)-\lambda(H^0_{\m}(S/JS))\\
&=& \lambda(S/x_{d-1}S: IS+ JS: \m^{\infty})- \lambda[(x_{d-1}S: IS+JS)\cap (JS: \m^{\infty})/JS]\\
&= & \lambda(S/x_{d-1}S: IS+ JS: \m^{\infty})\\
&=& \lambda(R/J_{d-1}:I+(J_{d-2}:I +I):\m^{\infty}).
  \end{eqnarray*}

  Finally for $n\geq 1$,
  \begin{eqnarray*}
  && \lambda[(x_{d-1}S)\cap I^nS/(x_{d-1}S)\cap I^{n+1}S+x_{d-1} I^{n-1}S]\\
&&-\lambda[I^{n+1}S:_{(x_{d-1}S:IS)\cap I^{n}S} \m^{\infty}/(x_{d-1}S:IS)\cap I^{n+1}S+x_{d-1}(I^{n}S:_{I^{n-1}S}\m^{\infty})]\\
&=& \lambda[(x_{d-1}S)\cap I^nS/(x_{d-1}S)\cap I^{n+1}S+x_{d-1} I^{n-1}S]\\
&&-\lambda[I^{n+1}S:_{(x_{d-1}S)\cap I^{n}S} \m^{\infty}/(x_{d-1}S)\cap I^{n+1}S+x_{d-1}(I^{n}S:_{I^{n-1}S}\m^{\infty})]\\
&\geq & 0, 
 \end{eqnarray*}
 since there is a map
 $$
 I^{n+1}S:_{(x_{d-1}S)\cap I^{n}S} \m^{\infty} \rightarrow (x_{d-1}S)\cap I^nS/(x_{d-1}S)\cap I^{n+1}S+x_{d-1} I^{n-1}S
 $$
 with kernel
\begin{eqnarray*}
&&[I^{n+1}S:_{(x_{d-1}S)\cap I^{n}S} \m^{\infty}]\cap [(x_{d-1}S)\cap I^{n+1}S+x_{d-1} I^{n-1}S]\\
&=& (x_{d-1}S)\cap I^{n+1}S+[I^{n+1}S:_{(x_{d-1}S)\cap I^{n}S} \m^{\infty}]\cap x_{d-1} I^{n-1}S\\
&=& (x_{d-1}S)\cap I^{n+1}S+[x_{d-1}I^{n}S:_{(x_{d-1}S)\cap I^{n}S} \m^{\infty}]\cap x_{d-1} I^{n-1}S\\
&=& (x_{d-1}S)\cap I^{n+1}S+x_{d-1}(I^{n}S:_{I^{n-1}S}\m^{\infty}),
\end{eqnarray*}
where the second equality holds because $\lambda(I^{n+1}S/x_{d-1}I^{n}S)<\infty$. 
\QED

\bigskip

The following theorem shows that the ideal $I$ enjoys nice properties when equality holds in the above inequality.
It generalizes the classical result of \cite{Huneke} and \cite{O}.

  \begin{Theorem}
  Assume $R$ is Cohen-Macaulay. let   $I$ be  an $R$-ideal with $\ell(I)=d$,  $G_d$ condition,
  $AN^-_{d-2}$ and ${\rm depth}(R/I)\geq {\rm min}\{1, {\rm dim}\,R/I\}$.
  Then for a general minimal reduction $J=(x_1, \ldots, x_d)$   of $I$, one has that
$j_1(I)=\lambda(I/J)+\lambda[R/J_{d-1}:I+(J_{d-2}:I+I):\m^{\infty}]$  if and only if $r(I)\leq 1$.
  In this case,  the associated graded ring $G$ is Cohen-Macaulay.
 \end{Theorem}
 
 \demo
 By the proof of Theorem \ref{North}, if $j_1(I)=\lambda(I/J)+\lambda[R/J_{d-1}:I+(J_{d-2}:I+I):\m^{\infty}]$ then
 for every $n\geq 1$, the length $\lambda[I^{n+1}S/JI^{n}S+ (x_{d-1}S:IS)\cap I^{n+1}S]=0$.
 Hence
 $$
 I^2\subseteq JI+ (J_{d-1}:I)\cap I^2=JI
 $$
 since $(J_{d-1}:I)\cap I^2=J_{d-1}I$ by \cite[Lemma 3.2]{PX}. Now the desired result follows from 
 \cite[Theorem~3.1]{JU}.
 \QED

\bigskip
   \begin{Corollary}
   Assume $R$ is Cohen-Macaulay. let  $I$ be  an $R$-ideal with $\ell(I)=d$,  $G_d$ condition
  and weakly $(d-2)$ residually $(S_2)$. Then for a general minimal reduction  $J=(x_1, \ldots, x_d)$ of $I$, one has

   (a)  $j_1(I)\geq 0$.

   (b) $j_1(I)= \lambda[R/J_{d-1}:I+(J_{d-2}:I+I):\m^{\infty}]$ implies that $I=J$ is a minimal reduction.

   (c) Assume $R$ is excellent. Then $j_1(I)= \lambda(I/J)$ implies that $I$ is $\m$-primary.

   (d) Assume $R$ is excellent. Then  $j_1(I)=0$ if and only if $I$ is a complete intersection.

   \end{Corollary}
   
   \demo
   (a) and (b) are clear. Assume (c). Then $\lambda[R/J_{d-1}:I+(J_{d-2}:I+I):\m^{\infty}]=0$, which implies
   $J_{d-1}:I+(J_{d-2}:I+I):\m^{\infty}=R$. Since $\ell(I)=d$, one has $J_{d-1}:I\neq R$.
   Hence $(J_{d-2}:I+I):\m^{\infty}=R$, i.e., ${\rm ht}(J_{d-2}:I+I)=d$. Since $R$ is excellent
  by \cite{NPUX}, \, ${\rm ht}(J_{d-2}:I+I)={\rm max}\{{\rm ht}\,I, d-1\}=d$, which yields ${\rm ht}\,I=d$, i.e., $I$ is  $\m$-primary. 
   The assertion (d) follows by (b) and (c).
   \QED
   
   \bigskip

We remark that (a) and (d) recover the work on the positivity of $j_1(I)$ by G. Colom$\acute{{\rm e}}$-Nin, C. Polini, B.  Ulrich and Y. Xie \cite{NPUX}.

We will finish the paper by an example from \cite{NPUX} which shows that if residual properties do not satisfy then the generalized Northcott's inequality fails to hold.

\begin{Example}\label{family}
{\rm Let $R=k\llbracket{x,y}\rrbracket/(x^3-x^2y)$ and $J=(xy^t)$ for any $t\ge 0$.  Notice that $R$ is an
 one-dimensional  Cohen-Macaulay local ring
and  $\ell(J)=1$. However,  $J$ does not satisfy $G_1$. By Macaulay2 \cite{M2},  one sees that $j_0(J)=t+1$, $j_1(J)=2-t$, which is strictly less than $0$ if $t>2$.}
\end{Example}

\end{document}